\documentclass[twoside,psfig]{article}

\usepackage{amsfonts,amssymb}

\oddsidemargin=\evensidemargin
\catcode`\@=11
\long\def\@makefntext#1{
\protect\noindent \hbox to 3.2pt {\hskip-.9pt  
$^{{\eightrm\@thefnmark}}$\hfil}#1\hfill}		

\def\ps@myheadings{\let\@mkboth\@gobbletwo		
\def\@oddhead{\hbox{}
\rightmark\hfil\eightrm\thepage}   
\def\@oddfoot{}\def\@evenhead{\eightrm\thepage\hfil
\leftmark\hbox{}}\def\@evenfoot{}
\def\sectionmark##1{}\def\subsectionmark##1{}}

\catcode`@=11		      		     
\def\ps@plain{\let\@mkboth\@gobbletwo
     \def\@oddhead{}\def\@oddfoot{\eightrm\hfil\thepage
     \hfil}\def\@evenhead{}\let\@evenfoot\@oddfoot}

\newcounter{sectionc}\newcounter{subsectionc}\newcounter{subsubsectionc}
\renewcommand{\section}[1] {\vspace{12pt}\addtocounter{sectionc}{1} 
\setcounter{subsectionc}{0}\setcounter{subsubsectionc}{0}\noindent 
	{\tenbf\thesectionc. #1}\par\vspace{5pt}}
\renewcommand{\subsection}[1] {\vspace{12pt}\addtocounter{subsectionc}{1} 
	\setcounter{subsubsectionc}{0}\noindent 
	{\bf\thesectionc.\thesubsectionc. 
	{\kern1pt \bfit #1}}\par\vspace{5pt}}
\renewcommand{\subsubsection}[1] {\vspace{12pt}
	\addtocounter{subsubsectionc}{1}
	\noindent
	{\tenrm\thesectionc.\thesubsectionc.\thesubsubsectionc.	{\kern1pt 
	\it #1}}\par\vspace{5pt}}
\newcommand{\nonumsection}[1] {\vspace{12pt}\noindent{\tenbf #1}
	\par\vspace{5pt}}

\newcommand{\Z}{\mathbb{Z}}

\newcounter{appendixc}
\newcounter{subappendixc}[appendixc]
\newcounter{subsubappendixc}[subappendixc]

\renewcommand{\appendix}[1] {\vspace{12pt}	
	\refstepcounter{appendixc}		
	\setcounter{figure}{0}
	\setcounter{table}{0}
	\setcounter{lemma}{0}
	\setcounter{theorem}{0}
	\setcounter{corollary}{0}
	\setcounter{definition}{0}
	\setcounter{equation}{0}
	\renewcommand{\thefigure}{\Alph{appendixc}.\arabic{figure}}
	\renewcommand{\thetable}{\Alph{appendixc}.\arabic{table}}
	\renewcommand{\theappendixc}{\Alph{appendixc}}
	\renewcommand{\thelemma}{\Alph{appendixc}.\arabic{lemma}}
	\renewcommand{\thetheorem}{\Alph{appendixc}.\arabic{theorem}}
	\renewcommand{\thedefinition}{\Alph{appendixc}.\arabic{definition}}
	\renewcommand{\thecorollary}{\Alph{appendixc}.\arabic{corollary}}
	\renewcommand{\theequation}{\Alph{appendixc}.\arabic{equation}}
	\noindent{\tenbf Appendix \theappendixc #1}\par\vspace{5pt}}

\newcommand{\textlineskip}{\baselineskip=13pt}
\newcommand{\smalllineskip}{\baselineskip=10pt}

\newcommand{\copyrightheading}[1]
	{\vspace*{-2.5cm}\smalllineskip{\flushleft
	{\footnotesize Journal of Knot Theory and Its Ramifications #1}\\
   	{\footnotesize \copyright\kern2pt World Scientific 
         Publishing Company}\\
         }}

\def\abstracts#1#2#3#4{{
	\centering{\begin{minipage}{4.5in}\footnotesize\baselineskip=10pt
	\centerline{ABSTRACT} 
	\parindent=15pt #1\par 
	\parindent=15pt #2\par
	\parindent=15pt #3\par
	\parindent=15pt #4\par
	\end{minipage}}\par}} 

\def\keywords#1{{ 
	\centering{\begin{minipage}{4.5in}\footnotesize\baselineskip=10pt
	{\footnotesize\it Keywords}\/: #1
	\end{minipage}}\par}}

\renewenvironment{thebibliography}[1]
	{\frenchspacing
	 \ninerm\baselineskip=11pt
	 \begin{list}{[\arabic{enumi}]}
	{\usecounter{enumi}\setlength{\parsep}{0pt}
	 \setlength{\leftmargin 13.7pt}{\rightmargin 0pt} 
	 \setlength{\itemsep}{0pt} \settowidth
	{\labelwidth}{[#1]}\sloppy}}{\end{list}}

\newcounter{itemlistc}
\newcounter{romanlistc}
\newcounter{alphlistc}
\newcounter{arabiclistc}

\newcommand{\fcaption}[1]{
        \refstepcounter{figure}
        \setbox\@tempboxa = \hbox{\footnotesize Fig.~\thefigure. #1}
        \ifdim \wd\@tempboxa > 5in
           {\begin{center}
        \parbox{5in}{\footnotesize\smalllineskip Fig.~\thefigure. #1}
            \end{center}}
        \else
             {\begin{center}
             {\footnotesize Fig.~\thefigure. #1}
              \end{center}}
        \fi}

\newcommand{\tcaption}[1]{
        \refstepcounter{table}
        \setbox\@tempboxa = \hbox{\footnotesize Table~\thetable. #1}
        \ifdim \wd\@tempboxa > 5in
           {\begin{center}
        \parbox{5in}{\footnotesize\smalllineskip Table~\thetable. #1}
            \end{center}}
        \else
             {\begin{center}
             {\footnotesize Table~\thetable. #1}
              \end{center}}
        \fi}

\def\pmb#1{\setbox0=\hbox{#1}
	\kern-.025em\copy0\kern-\wd0
	\kern.05em\copy0\kern-\wd0
	\kern-.025em\raise.0433em\box0}

\def\fnt#1#2{\footnotetext{\kern-.3em
	{$^{\mbox{\scriptsize #1}}$}{#2}}}

\def\runninghead#1#2{\pagestyle{myheadings}
\markboth{{\protect\footnotesize\it{\quad #1}}\hfill}
{\hfill{\protect\footnotesize\it{#2\quad}}}}

\font\tenrm=cmr10
 
\font\tenbf=cmbx10
\font\bfit=cmbxti10 at 10pt
\font\ninerm=cmr9

\font\eightrm=cmr8

\def\@begintheorem#1#2{\trivlist	
	\item[\hskip\labelsep{\bf #1\ #2.}]} 
\def\@opargbegintheorem#1#2#3{\trivlist
	\item[\hskip\labelsep{\bf #1\ #2\ (#3).}]}

    	{\setcounter{itemlistc}{0}		
	 \begin{list}{$\bullet$}		
	{\usecounter{itemlistc}			
	 \leftmargin10pt	       
	 \setlength{\parsep}{0pt}
	 \setlength{\itemsep}{0pt}     
	}}{\end{list}}

	{\setcounter{romanlistc}{0}		
	 \begin{list}{$($\roman{romanlistc}$)$}	
	{\usecounter{romanlistc}		
	 \leftmargin18pt 
	 \setlength{\parsep}{0pt}
	 \setlength{\itemsep}{0pt}	
	 \settowidth{\labelwidth}{#1}                          
	}}{\end{list}}

	{\setcounter{enumii}{0}			
	 \begin{list}{$($\alph{enumii}$)$}	
	{\usecounter{enumii}			
	 \leftmargin18pt		
	 \setlength{\parsep}{0pt}
	 \setlength{\itemsep}{0pt}	
	 \settowidth{\labelwidth}{#1}                          
	}}{\end{list}}

\def\qed{\hbox{${\vcenter{\vbox{			
   \hrule height 0.4pt\hbox{\vrule width 0.4pt height 6pt
   \kern5pt\vrule width 0.4pt}\hrule height 0.4pt}}}$}}

\def\theequation{\thesectionc.\arabic{equation}}  
\begin{document}

\setlength{\textheight}{7.7truein}  

\runninghead{J.S. Carter and M. Saito}{YBE and Fox calculus} 
\normalsize\textlineskip
\thispagestyle{empty}
\setcounter{page}{1}


\vspace*{0.88truein}

\centerline{\bf SET-THEORETIC YANG-BAXTER SOLUTIONS} 
\vspace{1mm}
\centerline{\bf VIA FOX CALCULUS}
\baselineskip=13pt
\vspace*{0.37truein}


\centerline{\footnotesize J. SCOTT CARTER}
\baselineskip=12pt
\centerline{\footnotesize\it Department of Mathematics, 
University of South Alabama,} 
\baselineskip=10pt
\centerline{\footnotesize\it Mobile, AL $36688$,    U.S.A.}
\baselineskip=10pt
\centerline{\footnotesize\tt carter@jaguar1.usouthal.edu}

\vspace*{10pt}
\centerline{\footnotesize MASAHICO SAITO}
\baselineskip=12pt
\centerline{\footnotesize\it Department of Mathematics, 
University of South Florida,} 
\baselineskip=10pt
\centerline{\footnotesize\it Tampa, FL $33620$, U.S.A.}
\baselineskip=10pt
\centerline{\footnotesize\tt saito@math.usf.edu}

\vspace*{20pt}

\centerline{\it Dedicated to  Professor Louis H. Kauffman for his $60$th birthday}



\vspace*{1in}

\abstracts{ We construct solutions to the set-theoretic Yang-Baxter equation
 using braid group representations in  free group automorphisms
 and their Fox  differentials. 
 The method resembles the extensions of groups and quandles.}{}{}{}

\vspace*{10pt}


\vspace*{1pt}\textlineskip



\input{epsf.sty}


\vspace{10mm}

For a set $X$, a mapping $R: X \times X \rightarrow  X \times X$
is called a solution to the set-theoretic Yang-Baxter equation 
(SYBE for short), if it 
satisfies the relation 
$$(R \times 1)(1 \times R)(R \times 1)=(1 \times R)(R \times 1)(1 \times R) $$
where $1$ denotes the identity map. 
Often $R$ is required to be invertible, which we do not impose in this paper,
and we concentrate on the above equation.
The set theoretic  Yang-Baxter equations are
 studied in detail in the papers, for example,  \cite{CESbiq,ESG,ESS,LYZ,Solo}. 
They have been subjects of active investigations, not only from 
quantum algebra points of view, but also for applications to knot theory.
In particular, algebraic systems called racks~\cite{FR,K&P}, 
quandles~\cite{Joyce,Matveev}, biracks, biquandles~\cite{BF,FJKW}
have been actively studied  in recent years in knot theory
(see also \cite{CJKLS,CKS,SilWil}, for example).  
It is of interest for applications to knots to have ample examples 
of such algebraic structures at hand to compute and use for
knot invariants. The purpose of this note is to construct examples
using representations of braid groups to the automorphisms of  
free groups~\cite{Wada}
 and the Fox free  differential calculus~\cite{Fox}.
 The construction resembles constructions of group 
 and quandle~\cite{CENS} extensions by $2$-cocycles.

Let $F_n=
\langle x_1, \ldots, x_n \rangle $ 
be the free group of rank $n$. 
A (Fox) derivative~\cite{Fox} is a $1$-cocycle, 
i.e., a homomorphism $d: \Z F_n \rightarrow  \Z F_n$
satisfying the $1$-cocycle condition (the Leibniz rule):
$d(uv)=d(u) + u\ d(v)$ for any $u, v \in F_n$. 
Fox showed that any derivative is uniquely written 
as a linear combination of 
$$ \frac{\partial}{\partial x_i}= \partial_{x_i} = \partial_i$$
defined (through linear extension) by
 $ \partial_{x_i} (x_j) = \delta(x_i, x_j)$ for $i, j=1, \ldots, n$,
 where $\delta$ denotes the Kronecker's delta.

The chain rule for Fox calculus~\cite{Fox} is 
writen as follows. Let
 $$\lambda: Y=\langle y_1, \ldots, y_{\ell} \rangle \rightarrow
  X=\langle x_1 \ldots, x_n \rangle $$
  be a free group
homomorphism.
Then for any $f \in \Z Y$, 
$$ \frac{\partial f^{\lambda} }{\partial x_j} 
= \sum _k \left(  \frac{\partial f }{\partial y_k} \right)^{\lambda} 
  \frac{\partial y_k ^{\lambda}  }{\partial x_j} ,$$
where $f^\lambda $ is defined as $ \hat{ \lambda} (f) $ by the homomorphism
$\hat{\lambda}: \Z Y \rightarrow \Z X$ induced from $\lambda$.

Wada~\cite{Wada}  
considered representations of braid groups to the automorphism 
groups of  free groups  of the following type:
For a standard  generator 
$\sigma_i$ of the $n$-string braid group $B_n$, an isomorphism of 
$\langle x_1, \ldots, x_n \rangle $ 
written as 
$$(x_i, x_{i+1}) \mapsto (u(x_i, x_{i+1}), v(x_i, x_{i+1})), \quad (i=1, \ldots, n-1) $$
is assigned, where $u,v$ are finite words in $x_i$, $i=1, \ldots, n$.
Such representations are also studied independently by A. J. Kelly~\cite{Kelly}
as mentioned in the review article by Przytycki~\cite{Pr} of Wada's paper.

\begin{figure}[htb]
\begin{center}
\mbox{
\epsfxsize=4in
\epsfbox{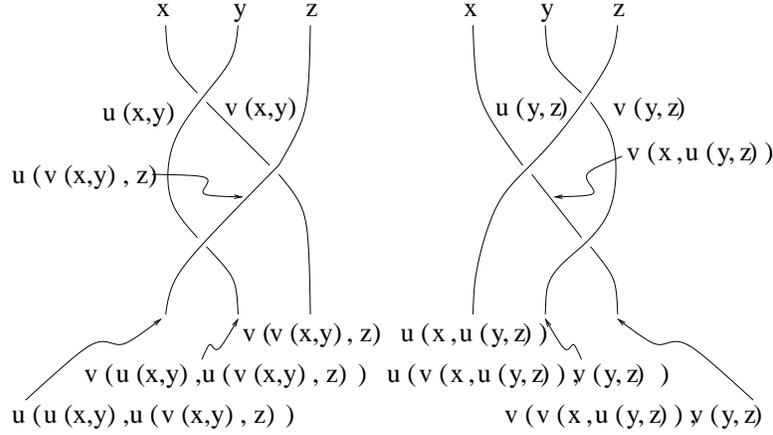} 
}
\end{center}
\caption{Braid relation and a biquandle axiom}
\label{biq} 
\end{figure}

The equalities for the braid relation 
$\sigma_i \sigma_{i+1} \sigma_i=\sigma_{i+1} \sigma_i \sigma_{i+1}$
 are written as 
$$\begin{array}{lccc}
(1) \quad  \quad \quad \quad \quad & u (u (x,y), u ( v (x, y), z) ) &=& u (x, u(y, z)) \\[2pt]
(2)  \quad \quad \quad \quad \quad  & v (u (x,y), u ( v (x, y), z) ) &=& u ( v (x, u (y, z) ), v (y,z) ) \\[2pt]
(3)  \quad \quad \quad  \quad \quad & v (v (x,y), z) &=& v ( v (x, u (y, z) ), v (y,z) )
\end{array}
$$
where $(x,y,z)=(x_i, x_{i+1}, x_{i+2})$ and $\sigma_i$ denotes a standard braid generator.
This computation is represented by diagrams as depcted in Fig.~\ref{biq}. 
If $u$ and $ v \in F(x,y)$ satisfy the above equalities (1), (2) and (3), 
then any group $G$ becomes a Yang-Baxter set 
with the operation $R: G^2 \rightarrow G^2$
defined by  $R(x,y)=(u(x, y), v(x, y))$ for $x, y \in G$.

Wada classified such representations for words 
of up to  word lengths $10$ for $u$ and $v$.
His list  of such functions consists of 
$(u, v)=(x,y)$, $(y^{-1}, x)$, $(y^{-1}, x^{-1})$, 
$(y, y^m x y^{-m})$, $(y, yx^{-1}y)$, 
$(y^{-1}, yxy)$, $(x^{-1} y^{-1} x, y^2 x)$.

Define
$$
\begin{array}{llllll}
u_1 (x,y)& =& \frac{\partial }{\partial x} u (x,y), \;
&u_2 (x,y)&=&\frac{\partial }{\partial y} u (x,y), \\[2mm]
v_1 (x,y)&=&\frac{\partial }{\partial x} v (x,y), \;
&v_2 (x,y)&=&\frac{\partial }{\partial y} v (x,y). 
\end{array}$$
These are elements of $\Z F_2$. 

\bigskip

{\bf Lemma 1.\/} 
{\it
Suppose $u, v \in F(x,y)$ satisfy the equations $(1)$, $(2)$, and $(3)$. 

\noindent
Then  $u_i, v_i$, $i=1,2$,  satisfy the following equalities in $\Z F_3$.
}

\medskip

\centerline{
$u_1(u (x,y), u ( v (x, y), z) ) u_1(x,y) 
+ u_2(u (x,y), u ( v (x, y), z) )  u_1(v(x,y),z) v_1(x,y)   $}

\centerline{
$=u_1(x, u(y,z)) , $}

\medskip

\centerline{
$
u_1(u (x,y), u ( v (x, y), z) ) u_2(x,y) 
+ u_2(u (x,y), u ( v (x, y), z) )  u_1(v(x,y),z) v_2(x,y)  $}
\centerline{
$=  u_2(x, u(y,z)) u_1(y,z) , $}

\medskip

\centerline{
$u_2(u (x,y), u ( v (x, y), z) )  u_2(v(x,y),z)  
=  u_2(x, u(y,z)) u_2(y,z),  $}

\bigskip 

\centerline{
$v_1(u (x,y), u ( v (x, y), z) ) u_1(x,y) 
+ v_2(u (x,y), u ( v (x, y), z) )  u_1(v(x,y),z) v_1(x,y)  $}
\centerline{
$= u_1(v (x, u (y, z) ), v (y,z) )  v_1(x, u(y,z)), $}

\medskip

\centerline{
$v_1(u (x,y), u ( v (x, y), z) ) u_2(x,y) 
+ v_2(u (x,y), u ( v (x, y), z) )  u_1(v(x,y),z) v_2(x,y)  $}
\centerline{
$ =   u_1(v (x, u (y, z) ), v (y,z) )  v_2(x, u(y,z)) u_1(y,z)
+ u_2(v (x, u (y, z) ), v (y,z) )  v_1(y,z) ,$}

\medskip

\centerline{
$ v_2(u (x,y), u ( v (x, y), z) ) u_2(v(x,y),z) = $}
\centerline{
$u_1( v (x, u (y, z) ), v (y,z) )v_2(x, u (y, z))u_2(y,z)  + 
u_2( v (x, u (y, z) ), v (y,z) )v_2(y,z) ,$}

\bigskip 

\centerline{
$v_1(v (x,y), z) v_1(x,y) = 
v_1( v (x, u (y, z) ), v (y,z) ) v_1(x, u (y, z) ), $}

\medskip

\centerline{
$ v_1(v (x,y), z) v_2(x,y) $}
\centerline{
$= v_1( v (x, u (y, z) ), v (y,z) ) v_2(x, u (y, z) ) u_1(y,z)
+ v_2( v (x, u (y, z) ), v (y,z) ) v_1(y,z) ,$}

\medskip

\centerline{
$  v_2(v (x,y), z) = $}
\centerline{
$ v_1( v (x, u (y, z) ), v (y,z) ) v_2(x, u (y, z) ) u_2(y,z)
+  v_2( v (x, u (y, z) ), v (y,z) ) v_2(y,z) . $}

\medskip

{\bf Proof.\/}
First consider the right-hand side $ u (x, u(y,z)) $ of
the equation (1). 
  Let  $\lambda : F( a, b ) \rightarrow 
 F (x,y,z )$
be defined by 
$\lambda (a)=x$ and  $\lambda(b)=u(y,z)$ that corresponds to 
the composition  $u (x, u(y, z)) $.
Applying the chain rule, we obtain  
\begin{eqnarray*}  \frac{\partial }{\partial y} u (x, u(y, z))& =&
\left(  \frac{\partial u(a,b) }{\partial a} \right)^{\lambda}
  \left(  \frac{\partial x }{\partial y} \right)
+  \left(  \frac{\partial u(a,b) }{\partial b} \right)^{\lambda}
 \left(  \frac{\partial u(y,z) }{\partial y} \right) \\[2mm]
 & =& u_2 (x, u(y, z) ) u_1(y,z) .
 \end{eqnarray*}
 By similar calculations for the equations (1), (2), (3) for
 derivatives with respect to $x$, $y$ and $z$, respectively, we obtain the above nine 
 equalities.
$\hfill$ $\Box$

\bigskip

Let $u, v \in F(x,y)$ be solutions to the equations (1), (2) and (3).
Define  maps $\hat{u}, \hat{v}: ( G \times V )^2 \rightarrow G \times V$
by 
\begin{eqnarray*}
\hat{u}( (x, a) , (y, b)) &=& (u(x,y), u_1(x,y) \ a + u_2(x,y)\  b ),  \\
\hat{v}( (x, a) , (y, b)) &=& (v(x,y), v_1(x,y)\  a + v_2(x,y)\  b )
\end{eqnarray*}
for $x,y \in  G$ and $a, b \in V$, where $V$ is a $G$-module. 

\bigskip

{\bf Theorem 2.\/} 
{\it Let $u, v \in F(x,y)$ be solutions to the equations $(1)$, $(2)$ and $(3)$.
Then the  map $R=(\hat{u}, \hat{v}) :  ( G \times V )^2 \rightarrow  ( G \times V )^2$
defined by 
$$R( (x, a) , (y, b) )=(\; \hat{u} ((x, a) , (y, b)), \; \hat{v}( (x, a) , (y, b)) \; ) $$ 
is a SYBE solution.
}

\medskip

{\bf Proof.\/}
We compute 

\medskip

\centerline{
$(R \times 1)(1 \times R)(R \times 1) ( (x,a), (y,b), (z,c) ) = $}

\centerline{ $( \; ( u (u (x,y), u ( v (x, y), z) ) , A), \
( v (u (x,y), u ( v (x, y), z) ) , B), \
( v (v (x,y), z) , C) \; )$ , }

\medskip

\centerline{$ (1 \times R)(R \times 1)(1 \times R) ( (x,a), (y,b), (z,c) ) = $}
\centerline{$ ( \; ( u (x, u(y, z)) , A'),\
 ( u ( v (x, u (y, z) ), v (y,z) ) , B'), \
 ( v ( v (x, u (y, z) ), v (y,z) ), C') \; ) 
 $}
 
 \medskip
 
 \noindent
where the second factors are as follows.
\begin{eqnarray*}
A&= & u_1  (u (x,y), u ( v (x, y), z) ) [ u_1 (x,y) \ a +  u_2(x,y)\ b\  ] \\
& & + u_2 (u (x,y),  u( v (x, y), z) ) \times \\
& & [ u_1 (v(x,y),z) (v_1(x,y) \ a + v_2 (x,y) \ b\  )
 + u_2 (v(x,y),z) \ c\  ]  , \\
 B &=& v_1  (u (x,y), u ( v (x, y), z) ) [ u_1 (x,y)\ a +  u_2(x,y)\ b\  ] \\
& & + v_2 (u (x,y),  u( v (x, y), z) ) \times \\
& & [ u_1 (v(x,y),z) (v_1(x,y) \ a + v_2 (x,y) \ b\  )
 + u_2 (v(x,y),z) \ c\  ] ,  \\
 C &=& v_1(v(x,y),z) [v_1(x,y)\ a  +  v_2(x,y)\ b\  ] + v_2(v(x,y),z)\ c ,
\end{eqnarray*}
\begin{eqnarray*}
A'&= & u_1  (x, u (  y, z) )  \ a +  u_2(x,u(y,z)) [
(u_1(y,z)\ b + u_2(y,z)\ c\ )  ] ,\\
B' &=& u_1(v (x, u(y,z)), v(y,z)) \times \\
& & [ v_1 (x,u(y,z)) \ a +  v_2(x,u(y,z))
(u_1(y,z)\ b + u_2(y,z)\ c\ ) ] ,\\
& &+ u_2(v (x, u(y,z)), v(y,z))[v_1(y,z)\ b + v_2(y,z)\ c \ ] \\
C'&=& v_1(v (x, u(y,z)), v(y,z)) \times \\
& & [ v_1 (x,u(y,z)) \ a +  v_2(x,u(y,z))
(u_1(y,z)\ b + u_2(y,z)\ c\ ) ] \\
& &+  v_2(v (x, u(y,z)), v(y,z))[v_1(y,z)\ b + v_2(y,z)\ c \ ] .
\end{eqnarray*}
By grouping coefficients for $a,b,c$ in the equations 
$A=A'$, $B=B'$, and $C=C'$, we obtain the result by Lemma~1.
$\hfill$ $\Box$

\bigskip

\begin{figure}[htb]
\begin{center}
\mbox{
\epsfxsize=3in
\epsfbox{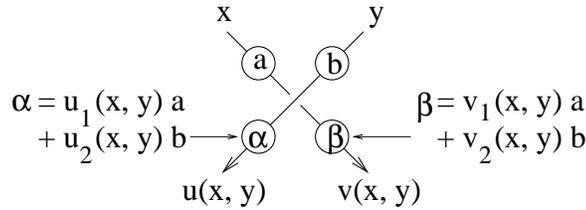} 
}
\end{center}
\caption{vectors assigned to colored arcs}
\label{YBEcrossing} 
\end{figure}

The above calculations can be visualized directly by diagrams of 
the Reidemeister type III move. For a given coloring by group elements 
at a crossing as in Fig.~\ref{YBEcrossing}, assign small beads at the top arcs
as depicted, representing elements 
($a, b$ in the figure) of a $G$-module $V$.
Below the crossing, elements of $V$ denoted by
$\alpha$ and $ \beta$ in the figure,  that are
 the images of 
linear maps $u_1(x,y)\ a+ u_2 (x,y)\ b$ and $v_1(x,y)\ a+ v_2 (x,y)\ b$,
are represented by another pair of beads as depicted.
Let elements $a, b, c \in V$ be assigned to the three arcs at the  top of Fig.~\ref{biq} 
from left to right, respectively. At the bottom, arcs receive three beads
representing the images of
 the composition of linear maps, corresponding to three crossings. 
 The outcomes are $A$, $B$, $C$ in the above calculation
for the left-hand side, and $A'$, $B'$, $C'$ for the right-hand side.

{\bf Example 3.\/} 
For the second  last example in Wada's list
$u=y^{-1} $ and $v=yxy$, 
we obtain 
$$ u_1=0, \quad 
u_2=-y^{-1}, \quad 
v_1=y, \quad
v_2=1+yx ,
$$
so that for any group $G$ and any $G$-module $V$,
$G \times V$ is a Yang-Baxter set by 
$$
R((x,a), (y,b)) 
= (\  (y^{-1} ,  (- y^{-1})b \ ),\ 
(y xy, y a + (1+yx)b \ ) \ ) .
$$
Similarly, for the last example in Wada's list
$u=x^{-1}y^{-1} x$ and $v=y^2 x$, 
we obtain 
\begin{eqnarray*}
\lefteqn{
R((x,a), (y,b)) } \\
&=& (\  (x^{-1}y^{-1} x, (- x^{-1}+ x^{-1}y^{-1})a + (-x^{-1} y^{-1})b \ ), \ 
(y^2 x, (y^2)a + (1+y)b \ ) \ ) .
\end{eqnarray*}

\bigskip

{\bf Remark 4.\/} 
For the $4$th and $5$th examples in the Wada's list, 
the operation $x*y=v(x,y)$ (where  $u(x,y)=y$ does not play an essential role) 
defines a rack, 
and the equalities in Lemma~1 
give rise to the condition of rack  algebras and modules defined 
in \cite{AG}; in their notation, 
$\eta_{x,y}=v_1(x,y)$ and $\tau_{x,y}=v_2(x,y)$. Hence in these cases, Theorem~2 
 can be used to produce examples of rack modules.
 
 On the other hand, the equalities in Lemma~1 
 can be regarded  as providing a natural definition of birack algebras and modules,
 from point of view of \cite{AG}. 
 Specifically, a birack algebra can be defined with generators 
$$\{ u_i(x,y), v_j(x,y) \ | \ i, j=1, 2, \ x, y \in X \} $$
 with relations stated in Lemma~1, 
 where $X$ is a birack defined in \cite{BF,FJKW}. 
 Furthermore, again from \cite{AG}, it is expected that 
 these could be used as twisted coefficients for generalizations 
 of the homology theory of the SYBE~\cite{CESbiq}.
 Applications to knot invariants are also expected.

\nonumsection{Acknowledgements}
The first and the second authors are
partially supported by 
NSF Grant DMS $\#0301095$ and 
 $\#0301089$, respectively.

\nonumsection{References}

\end{document}